\documentstyle[leqno]{book}
\newfont{\newit}{cmfi10 scaled 1100}
\newfont{\fnewit}{cmfi10 scaled 800}
\newfont{\cirilrm}{wncyr10 scaled 1000}
\newfont{\cirilbf}{wncyb10 scaled 1000}
\newfont{\cirilsf}{wncyss10 scaled 1000}
\newfont{\cirilit}{wncyi10 scaled 1000}
\newfont{\cirilsc}{wncysc10 scaled 1000}

\newcommand{\z}{\symbol{'31}}

\newcommand{\sh}{\symbol{'170}}

\newcommand{\ja}{\symbol{'37}}

\newcommand{\ij}{\symbol{'32}}

\newcommand{\ii}{\symbol{'171}}

\newcounter{supersection}[section]
\newtheorem{th}[supersection]{Theorem}
\newtheorem{lm}[supersection]{Lemma}

\newtheorem{co}[supersection]{Corollary}

\def\bibname{\textbf{REFERENCES}}
\def\thebibliography#1{\paragraph*{\uppercase{\bibname}}\list
{[\arabic{enumi}]}{\settowidth\labelwidth{[#1]}\leftmargin\labelwidth
\advance\leftmargin\labelsep\usecounter{enumi}}
\def\newblock{\hskip .11em plus .33em minus .07em}
\sloppy\clubpenalty4000\widowpenalty4000
\sfcode`\.=1000\relax}

\def\Dj{D{\hspace{-.75em}\raisebox{.3ex}{-}\hspace{.4em}}}
\def\stop{\mbox{\footnotesize {\vrule width 6pt height 6pt}}}

\begin{document}

\thispagestyle{plain}

\bigskip

\centerline{\large \bf SOME INEQUALITIES FOR}

\smallskip

\centerline{\large \bf ALTERNATING KUREPA'S FUNCTION}
\footnotetext{Research partially supported by the MNTRS,
Serbia \& Montenegro, Grant No. 1861.}

\vspace*{6.00 mm}

\centerline{\large \it Branko J. Male\v sevi\' c}

\vspace*{4.00 mm}

\begin{center}
\parbox{25.0cc}{\scriptsize \bf \boldmath
In this paper we consider alternating Kurepa's function $A(z)$
\cite{Petojevic_02}. We give some recurrent relations for
alternating Kurepa's function via appropriate sequences of rational
functions and gamma function. Also we give some inequalities for the
real part of alternating Kurepa's function $A(x)$ for values of
argument $x \!>\! -2$. The obtained results are analogous to results
from \cite{Malesevc_04}.}
\end{center}

\noindent
\section{\large \bf \boldmath \hspace*{-7.0 mm}
1. Alternating Kurepa's function $A(z)$}

\smallskip
\noindent
{\sc R. Guy} considered, in the book \cite{Guy_94} (p.$\,$100.), the function
of alternating left factorial $A(n)$ as an alternating sum of factorials
$A(n) = n! - (n-1)! + \ldots + (-1)^{n-1}1!$. Let us use the standard
notation:
\begin{equation}
\label{A_SUM_1}
A(n) = \displaystyle\sum\limits_{i=1}^{n}{(-1)^{n-i}i!}.
\end{equation}
Sum (\ref{A_SUM_1}) corresponds to the sequence $A005165$ in \cite{Sloane_03}.
An analytical extension of the function (\ref{A_SUM_1}) over the set of complex
numbers is determined by the integral:
\begin{equation}
\label{A_INT_1}
A(z)
=
\displaystyle\int\limits_{0}^{\infty}{
e^{-t} \displaystyle\frac{t^{z+1}-(-1)^{z}t}{t+1} \: dt},
\end{equation}
which converges for $\mbox{Re} \: z \!>\! 0$ \cite{Petojevic_02}. For function $A(z)$ we use
the term {\it alternating Kurepa's function}. It is easily verified that alternating
Kurepa's function is a solution of the functional equation:
\begin{equation}
\label{A_FE_1}
A(z) + A(z-1) = \Gamma(z+1).
\end{equation}
Let us observe that since $A(z-1) = \Gamma(z+1) - A(z)$, it is possible
to make the analytical continuation of alternating Kurepa's function $A(z)$
for \mbox{$\mbox{Re} \, z \leq 0$}. In that way, the alternating Kurepa's function $A(z)$
is  a meromorphic function with simple poles at $z = -n$~\mbox{$(n \!\geq\! 2)$}
~\cite{Petojevic_02}.

\smallskip
\noindent
Let us emphasize that in the following consideration, in the sections
{\bf 2.} and {\bf 3.}, it is sufficient to use only fact that function $A(z)$
is a solution of the functional equation~(\ref{A_FE_1}). In section {\bf 4.}
we give some inequalities for the real part of alternating Kurepa's function $A(x)$
for values of argument $x \!>\! -2$.

\break

\noindent
\section{\large \bf \boldmath \hspace*{-7.0 mm}
2. Representation of the alternating Kurepa's function \\
\hspace*{-1.0  mm} via sequences of polynomials and gamma function}

Let us introduce a sequences of polynomials:
\begin{equation}
\label{Ap_n_Def}
\mbox{\newit p}_{n}(z) = (z-n+1)\mbox{\newit p}_{n-1}(z) + (-1)^{n},
\end{equation}
with initial member $\mbox{\newit p}_{0}(z)=1$. Analogously to results
from \cite{Kurepa_73}, the following statements are true\footnotetext{
Letters {\fnewit p}, {\fnewit q}, {\fnewit r}, {\fnewit g} are printed
in the funny italic {\TeX } font.}:
\begin{lm}
For each $n \!\in\! \mbox{\bf N}$ and $z \!\in\! \mbox{\bf C}$ we have explicitly$:$
\begin{equation}
\mbox{\newit p}_{n}(z)
\!=\!
(-1)^{n} {\Big (} 1 + \displaystyle\sum_{j=0}^{n-1}{
\displaystyle\prod_{i=0}^{j}{(-1)^{j-1}(z-n+i+1)}} {\Big )}.
\end{equation}
\end{lm}
\begin{th}
For each $n \!\in\! \mbox{\bf N}$ and $z \!\in\! \mbox{\bf C} \backslash
\mbox{\big (}\mathop{\mbox{\bf Z}}^{-} \cup \, \{0,1,2,\ldots,n\!-\!2\} \mbox{\big )}$
is valid$:$
\begin{equation}
\label{A_P_Veza}
A(z) = (-1)^{n}A(z - n) + \mbox{\newit p}_{n-1}(z) \cdot \Gamma(z\!-\!n\!+\!2).
\end{equation}
\end{th}

\section{\large \bf \boldmath \hspace*{-7.0 mm}
3. Representation of the alternating Kurepa's function \\
\hspace*{-1.0 mm} via sequences of rational functions and gamma function}

Let us observe that on the basis of a functional equation for the gamma function
$\Gamma(z+1) = z\Gamma(z)$, it follows that the alternating Kurepa's function is
solution of the following functional equation:
\begin{equation}
\label{A_FE_2}
A(z+1) - zA(z) - (z+1)A(z-1) = 0.
\end{equation}
For $z \!\in\! \mbox{\bf C} \backslash \{-1\}$, based on (\ref{A_FE_2}), we have:
\begin{equation}
\label{A_FE_2_1}
A(z-1)
=
-
\displaystyle\frac{z}{z+1}A(z)
+
\displaystyle\frac{1}{z+1}A(z+1)
=
\mbox{\newit q}_{1}(z) A(z) - \mbox{\newit r}_{1}(z) A(z+1),
\end{equation}
for rational functions $\mbox{\newit q}_{1}(z) \!=\! -\mbox{\small $\displaystyle\frac{z}{z\!+\!1}$}$,
$\mbox{\newit r}_{1}(z) \!=\! -\mbox{\small $\displaystyle\frac{1}{z\!+\!1}$}$ over $\mbox{\bf C}
\backslash \{-1\}$. Next, for $z \in \mbox{\bf C} \backslash \{-1,0\}$, based on (\ref{A_FE_2}),
we obtain:
\begin{equation}
\begin{array}{rcl}
A(z-2)
&\!\!\!\!=\!\!\!\!&
\displaystyle\frac{1}{z}A(z)
\!-\!
\displaystyle\frac{z\!-\!1}{z}A(z-1)                                                  \\[2.0 ex]
&\!\!\!\!\mathop{=}\limits_{(\ref{A_FE_2_1})}\!\!\!\!&
\displaystyle\frac{1}{z}A(z)\!-\!\displaystyle\frac{z\!-\!1}{z}{\bigg (}
\!-\!\displaystyle\frac{z}{z\!+\!1}A(z)\!+\!\frac{1}{z\!+\!1}A(z+1){\bigg )}          \\[2.0 ex]
&\!\!\!\!=\!\!\!\!&
\displaystyle\frac{z^2\!+\!1}{z(z\!+\!1)}A(z)
\!-\!
\displaystyle\frac{z\!-\!1}{z(z\!+\!1)}A(z+1)
=
\mbox{\newit q}_{2}(z) A(z)\!-\!\mbox{\newit r}_{2}(z) A(z+1),
\end{array}
\end{equation}
for rational functions $\mbox{\newit q}_{2}(z) \!=\! \mbox{\small $\displaystyle\frac{
z^2\!+\!1}{z(z\!+\!1)}$}$, $\mbox{\newit r}_{2}(z) \!=\! \mbox{\small $\displaystyle\frac{
z\!-\!1}{z(z\!+\!1)}$}$ over $\mbox{\bf C} \backslash \{-1,0\}$. Thus, for values
$z \in \mbox{\bf C} \backslash \{-1,0,1,\ldots,n\!-\!2\}$, based on (\ref{A_FE_2}),
by mathematical induction it is true:

\smallskip \noindent
\begin{equation}
A(z-n)
=
\mbox{\newit q}_{n}(z) A(z) - \mbox{\newit r}_{n}(z) A(z+1),
\end{equation}
for rational functions $\mbox{\newit q}_{n}(z)$, $\mbox{\newit r}_{n}(z)$ over $\mbox{\bf C}
\backslash \{-1,0,1,\ldots,n\!-\!2\}$ which fulfill the same recurrent relations:
\begin{equation}
\label{Rec_q}
\mbox{\newit q}_{n}(z)
=
-\displaystyle\frac{z-n+1}{z-n+2} \, \mbox{\newit q}_{n-1}(z)
+\displaystyle\frac{1}{z-n+2} \, \mbox{\newit q}_{n-2}(z)
\end{equation}
and

\begin{equation}
\label{Rec_r}
\mbox{\newit r}_{n}(z)
=
-\displaystyle\frac{z-n+1}{z-n+2} \, \mbox{\newit r}_{n-1}(z)
+\displaystyle\frac{1}{z-n+2} \, \mbox{\newit r}_{n-2}(z),
\end{equation}
with different initial functions $\mbox{\newit q}_{1,2}(z)$ and $\mbox{\newit r}_{1,2}(z)$.

\medskip
\noindent
Based on the previous consideration we can conclude:

\begin{lm}
For each $n \!\in\! \mbox{\bf N}$ and $z \!\in\! \mbox{\bf C} \backslash \{-1,0,1,\ldots,n\!-\!2\}$
let the rational function $\mbox{\newit q}_{n}(z)$ be determined by the recurrent relation
{\rm (\ref{Rec_q})} with initial functions $\mbox{\newit q}_{1}(z) \!=\!
-\mbox{\small $\displaystyle\frac{z}{z+1}$}$ and $\mbox{\newit q}_{2}(z) \!=\!
\mbox{\small $\displaystyle\frac{z^2+1}{z(z+1)}$}$. Thus the sequences
$\mbox{\newit q}_{n}(z)$ has an explicit form$:$
\begin{equation}
\mbox{\newit q}_{n}(z)
=
(-1)^{n}
{\bigg (}
1
+
\displaystyle\sum\limits_{j=1}^{n}{
\displaystyle\prod\limits_{i=1}^{j}{
\displaystyle\frac{(-1)^j}{z+2-i}}}
{\bigg )}.
\end{equation}
\end{lm}

\begin{lm}
For each $n \!\in\! \mbox{\bf N}$ and $z \!\in\! \mbox{\bf C} \backslash \{-1,0,1,\ldots,n\!-\!2\}$
let the rational function $\mbox{\newit r}_{n}(z)$ be determined by the recurrent relation
{\rm (\ref{Rec_r})} with initial functions $\mbox{\newit r}_{1}(z) \!=\!
-\mbox{\small $\displaystyle\frac{1}{z+1}$}$ and $\mbox{\newit r}_{2}(z) \!=\!
\mbox{\small $\displaystyle\frac{z-1}{z(z+1)}$}$. Thus the sequences
$\mbox{\newit r}_{n}(z)$ has an explicit form$:$
\begin{equation}
\mbox{\newit r}_{n}(z)
=
(-1)^{n-1}
{\bigg (}
\displaystyle\sum\limits_{j=1}^{n}{
\displaystyle\prod\limits_{i=1}^{j}{
\displaystyle\frac{(-1)^j}{z+2-i}}}
{\bigg )}.
\end{equation}
\end{lm}

\begin{th}
\label{Th}
For each $n \!\in\! \mbox{\bf N}$ and $z \!\in\! \mbox{\bf C} \backslash \{-1,0,1,\ldots,n\!-\!2\}$
we have$:$
\begin{equation}
A(z) = (-1)^{n} {\Big (} A(z - n)
+
{\big (}(-1)^{n} - \mbox{\newit q}_{n}(z){\big )} \cdot \Gamma(z+2){\Big )}
\end{equation}
and

\begin{equation}
A(z) = (-1)^{n} {\Big (} A(z - n) + \mbox{\newit r}_{n}(z) \cdot \Gamma(z+2){\Big )}.
\end{equation}
\end{th}

\newpage

\section{\large \bf \boldmath \hspace*{-7.0 mm}
4. Some inequalities for the real part of alternating Kurepa's
   function}

In this section we consider alternating Kurepa's function $A(x)$, given by an integral
representation (\ref{A_INT_1}), for values of argument $x \!>\! -2$. The real and imaginary
parts of the function $A(x)$ are represented by:
\begin{equation}
\mbox{\rm Re}\:A(x) = \displaystyle\int\limits_{0}^{\infty}{
e^{-t}\displaystyle\frac{t^{x+1}- \cos (\pi x)\,t}{t+1}\:dt}
\end{equation}
and

\medskip
\noindent
\begin{equation}
\mbox{Im}\:A(x) = -\displaystyle\int\limits_{0}^{\infty}{
e^{-t}\displaystyle\frac{\sin (\pi x)\,t}{t+1}\:dt}.
\end{equation}
In this section we give some inequalities for the real part of alternating Kurepa's
function $A(x)$ for values of argument $x \!>\! -2$. The following statements are true:

\begin{lm}
\label{lm_21}
The function$:$
\begin{equation}
\label{lm_21_eq_1}
\beta(x) = \displaystyle\int\limits_{0}^{\infty}{e^{-t}\displaystyle\frac{t^{x+1}}{t+1}\:dt},
\end{equation}
over set $(-2,\infty)$ is positive, convex and fulfill an inequality$:$
\begin{equation}
\label{lm_21_eq_3}
\beta(x) \!\geq\! \beta(x_{0}) = 0.401\:855\:\ldots\:,
\end{equation}
with equality in the point $x_{0} = - 0.108\:057\:\ldots\:$.
\end{lm}

\noindent
{\bf Proof.} For positive function $\beta(x) \!\in\! C^{2}(-2,\infty)$, on the basis
of (\ref{lm_21_eq_1}), the condition of convexity $\beta^{''}(x) \!>\! 0$
is true. Next, based on (\ref{lm_21_eq_1}), we can conclude
\mbox{$\lim_{\varepsilon \rightarrow 0+}{\beta(-2+\varepsilon)}
= +\infty $}
and
\mbox{$\lim_{x \rightarrow +\infty}{\beta(x)}
= +\infty $}.
Therefore, we can conclude that exists exactly one minimum $x_{0} \!\in\! (-2,+\infty)$.
Using standard numerical methods it is easily determined
$x_{0} = - 0.108\:057\:\ldots$ and $\beta(x_{0}) = 0.401\:855\:\ldots\:$.~\stop

\begin{lm}
\label{lm_22}
The function$:$
\begin{equation}
\label{lm_22_eq_1}
\gamma(x)
=
\displaystyle\int\limits_{0}^{\infty}{e^{-t}\displaystyle\frac{\cos (\pi x) \, t}{t+1}\:dt},
\end{equation}
over set $(-2,\infty)$, is determined with$:$

\begin{equation}
\label{lm_22_eq_2}
\gamma(x)
=
(1+\mbox{\rm e}\,\mbox{\rm Ei}(-1)) \cdot \cos (\pi x)
=
0.403 \: 652 \: \ldots \: \cdot \cos (\pi x).
\end{equation}

\noindent
where
$\mbox{\rm Ei}(t) =
\displaystyle\!\int\limits_{-\infty}^{t}{\!\!\mbox{\small $\displaystyle\frac{e^u}{u}$}\;du}
\;\; (t<0)$ is function of exponential integral $(${\rm \cite{Rizik_71}, 8.211-1}$)$.
\end{lm}

\begin{lm}
\label{lm_23}
The function $\mbox{\rm Re}\:A(x)$, over set $(-2,\infty)$, is determined as difference$:$
\begin{equation}
\label{lm_23_eq_1}
\mbox{\rm Re}\:A(x)
=
\beta(x) - \gamma(x).
\end{equation}
and has two roots $x_{1} = -0.015\:401\:\ldots\:$ and $x_{2} = 0$.
The function $\mbox{\rm Re}\:A(x)$ is positive over set$:$

\smallskip \noindent
\begin{equation}
\label{lm_23_eq_2}
D_{1} = (-2,x_{1}) \cup (0,\infty)
\end{equation}
and negative over set$:$

\smallskip \noindent
\begin{equation}
\label{lm_23_eq_3}
D_{2} = (x_{1},0).
\end{equation}
\end{lm}

\noindent
{\bf Proof.} Let $\beta(x)$ be function from lemma \ref{lm_21} and
let $\gamma(x)$ be function from lemma \ref{lm_22}. For value \mbox{$x_2 = 0$}
it is true $\beta(x_2) = \gamma(x_2) = 0.403\,652\,\ldots$, ie. value $x_2 = 0$
is a root of function $\mbox{\rm Re}\:A(x)$. Let us prove that function
$\mbox{\rm Re}\:A(x)$ has exactly one root $x_{1} \in (x_0,x_2)$, where
$x_0 = -0.108\,057\,\ldots$ is value from lema \ref{lm_21}. It is true
$\beta(x_0) = 0.401\,855\,\ldots > 0.380\,061\,\ldots = \gamma(x_0)$.
Let us notice that $\beta(x)$ is convex and increasing function over set
$(x_{0},x_{2})$ and let us notice that $\gamma(x)$ is concave and increasing
function over same set $(x_{0},x_{2})$. Therefore, we can conclude that function
$\mbox{\rm Re} \: A(x)$ has exactly one root $x_{1} \!\in\! (x_{0},x_{2})$.
Using numerical methods we can determined $x_{1} = -0.015\:401\:\ldots\:$.
On the basis of the graphs of the functions $\beta(x)$ and $\gamma(x)$ we can conclude
that function $\mbox{\rm Re} \: A(x)$ has exactly two roots $x_{1}$ and $x_{2}$ over
set $(-2,\infty)$. Hence, the sets $D_{1}$ and $D_{2}$ are correctly determined.~\stop

\begin{lm}
\label{lm_24}
For $x \!\in\! (-1,1+x_1] \cup [1,\infty)$ it is true$:$
\begin{equation}
\label{lm_24_eq_1}
\Gamma(x+1) \geq \mbox{\rm Re}\:A(x),
\end{equation}
while the equality is true for $x\!=\!1\!+\!x_{1}$ or $x\!=\!1$.
\end{lm}

\noindent
{\bf Proof.} For $x \!>\! -1$ it is true:
\begin{equation}
\Gamma(x+1) \geq \mbox{Re}\:A(x) = \Gamma(x+1) - \mbox{\rm Re}\:A(x-1)
\:\Longleftrightarrow\:
\mbox{Re}\:A(x-1) \geq 0.
\end{equation}
Right side of the previous equivalence is true for $x\!-\!1 \!\in\!
(-2,x_1] \cup [0,\infty)$, ie. \mbox{$x \!\in\! (-1,1+x_1] \cup [1,\infty)$}.~\stop

\medskip
\noindent
In the following considerations let us denote
$\mbox{\bf E}_{a} = (a,a\!+\!2\!+\!x_1] \cup [a\!+\!2,\infty)$
for fixed $a \,\!\geq\! -1$.

\begin{co}
\label{Lm}
For fixed $k \!\in\! \mbox{\bf N}$ and values $x \in \mbox{\bf E}_{k}$ following
inequality is true$:$
\begin{equation}
\label{Ineq_GA1}
\displaystyle\frac{\mbox{\rm Re}\:A(x\!-\!k\!-\!1)}{\Gamma(x\!-\!k)} \leq 1,
\end{equation}
while the equality is true for $x = k\!+\!2\!+\!x_1$ or $x = k\!+\!2$.
\end{co}

\medskip
\noindent
In the next two proofs of theorems which follows
we use the auxiliary sequences of functions:
\begin{equation}
\quad\quad\quad\quad\quad\quad\quad\quad\quad\quad
\mbox{\newit g}_{k}(x) = \sum_{i=0}^{k-1}{(-1)^{k+i}\Gamma(x+1-i)}
\quad\quad\quad\quad\quad\quad\quad (k \in \mbox{\bf N}),
\end{equation}
for values $x \!>\! k\!-\!2$. Let us notice that for $x \!>\! k\!-\!2$ it is true:
\begin{equation}
\label{A_gr_veza}
\mbox{\newit g}_{k}(x)\!=\!\Gamma(x+2) \cdot \mbox{\newit r}_{k}(x).
\end{equation}
Then, the following statements are true:

\begin{th}
\label{A_Th1}
For fixed odd number $k = 2n\!+\!1\!\in\! \mbox{\bf N}$ and values
$x \!\geq\! k\!+\!1$ the following double inequality is true$:$
\begin{equation}
\label{Ineq_GA2}
\displaystyle\frac{\mbox{\newit p}_{k}(x)}{\mbox{\newit p}_{k}(x)+1}
\cdot {\big (}-\mbox{\newit r}_{k}(x){\big )}
\leq
\displaystyle\frac{\mbox{\rm Re}\:A(x)}{\Gamma(x+2)}
<
{\big (}-\mbox{\newit r}_{k}(x){\big )},
\end{equation}
while the equality is true for $x = k\!+\!1$.
\end{th}

\smallskip \noindent
{\bf Proof.} Based on lemma \ref{lm_23}, using theorem \ref{Th},
the following inequality is true:
\begin{equation}
\label{A_nepar_1}
\mbox{\rm Re}\:A(x) \leq -\mbox{\newit g}_{2n+1}(x),
\end{equation}
for values $x\!\in\!\mbox{\bf E}_{k-2}$. On the other hand, based on
(\ref{Ineq_GA1}), for values $x\!\in\!\mbox{\bf E}_{k-1}$ we can conclude:
\begin{equation}
\label{A_nepar_2}
\begin{array}{rcl}
\displaystyle\frac{\mbox{\rm Re}\:A(x)}{\mbox{\newit g}_{2n+1}(x)}
& \!\!\!=\!\!\! &
-1\!+\!\displaystyle\frac{\mbox{Re}\:A(x\!-\!2n\!-\!1)}{\mbox{\newit g}_{2n+1}(x)}
\:=\:
-1\!+\!\displaystyle\frac{\mbox{Re}\:A(x\!-\!2n\!-\!1)}{\Gamma(x\!-\!2n)(\mbox{\newit p}_{2n+1}(x)
                                                                                     \!+\!1)}\\[2.5 ex]
& \!\!\!=\!\!\! &
-1\!+\!\displaystyle\frac{\mbox{Re}\:A(x\!-\!2n\!-\!1)/\Gamma(x\!-\!2n)}{\mbox{\newit p}_{2n+1}(x)+1}
\:\leq\:
-\displaystyle\frac{\mbox{\newit p}_{2n+1}(x)}{\mbox{\newit p}_{2n+1}(x)+1}.
\end{array}
\end{equation}
From (\ref{A_nepar_1}) and (\ref{A_nepar_2}), using (\ref{A_gr_veza}), the
double inequality (\ref{Ineq_GA2}) follows for values $x \!\geq\! k \!+\! 1$. \stop

\begin{th}
\label{A_Th2}
For fixed even number $k = 2n \!\in\! N$ and values $x \!\geq\! k\!+\!1$ the
following double inequality is true$:$
\begin{equation}
\label{Ineq_GA3}
\mbox{\newit r}_{k}(x)
<
\displaystyle\frac{\mbox{\rm Re}\:A(x)}{\Gamma(x+2)}
\leq
\displaystyle\frac{\mbox{\newit p}_{k}(x)}{\mbox{\newit p}_{k}(x)-1}
\cdot \mbox{\newit r}_{k}(x),
\end{equation}
while the equality is true for $x = k\!+\!1$.
\end{th}

\smallskip \noindent
{\bf Proof.} Based on lemma \ref{lm_23}, using theorem \ref{Th},
the following inequality is true:
\begin{equation}
\label{A_par_1}
\mbox{\rm Re}\:A(x) \geq \mbox{\newit g}_{2n}(x),
\end{equation}
for values $x\!\in\!\mbox{\bf E}_{k-2}$. On the other hand, based on
(\ref{Ineq_GA1}), for values $x\!\in\!\mbox{\bf E}_{k-1}$ we can conclude:
\begin{equation}
\label{A_par_2}
\begin{array}{rcl}
\displaystyle\frac{\mbox{\rm Re}\:A(x)}{\mbox{\newit g}_{2n}(x)}
& \!\!\!=\!\!\! &
1\!+\!\displaystyle\frac{\mbox{Re}\:A(x\!-\!2n\!)}{\mbox{\newit g}_{2n}(x)}
\:=\:
1\!+\!\displaystyle\frac{\mbox{Re}\:A(x\!-\!2n\!)}{\Gamma(x\!-\!2n\!+\!1)(\mbox{\newit p}_{2n}(x)
                                                                                     \!-\!1)}\\[2.5 ex]
& \!\!\!=\!\!\! &
1\!+\!\displaystyle\frac{\mbox{Re}\:A(x\!-\!2n)/\Gamma(x\!-\!2n\!+\!1)}{\mbox{\newit p}_{2n}(x)\!-\!1}
\:\leq\:
\displaystyle\frac{\mbox{\newit p}_{2n}(x)}{\mbox{\newit p}_{2n}(x)-1}.
\end{array}
\end{equation}
From (\ref{A_par_1}) and (\ref{A_par_2}), using (\ref{A_gr_veza}), the
double inequality (\ref{Ineq_GA3}) follows for values $x \!\geq\! k \!+\! 1$. \stop

\begin{co}
\label{A_Co}
For fixed number $k \!\in\! \mbox{\bf N}$ and values
$x \!\geq\! k\!+\!1$ the following double inequality is true$:$
\begin{equation}
\label{Ineq_GA4}
\mbox{\newit r}_{k}(x)
<
(-1)^{k} \displaystyle\frac{\mbox{\rm Re}\:A(x)}{\Gamma(x+2)}
\leq
\displaystyle\frac{\mbox{\newit p}_{k}(x)}{\mbox{\newit p}_{k}(x)-(-1)^{k}}
\cdot \mbox{\newit r}_{k}(x),
\end{equation}
while the equality is true for $x = k\!+\!1$.
\end{co}

\begin{co}
On the basis of theorems {\rm \ref{A_Th1}} and {\rm \ref{A_Th2}} we can conclude$:$
\begin{equation}
\lim\limits_{x \rightarrow \infty}{\displaystyle\frac{\mbox{\rm Re}\:A(x)}{\Gamma(x+2)}}=0
\qquad\mbox{and}\qquad
\lim\limits_{x \rightarrow \infty}{\displaystyle\frac{\mbox{\rm Re}\:A(x)}{\Gamma(x+1)}}=1.
\end{equation}
\end{co}

\bigskip

\bigskip
{\small
\noindent University of Belgrade,
           \hfill (Received October 20, 2004)        \break
\noindent Faculty of Electrical Engineering,   \hfill\break
%           \hfill (Revised ?? ??, ????)             \break
\noindent P.O.Box 35-54, $11120$ Belgrade,     \hfill\break
\noindent Serbia \& Montenegro                 \hfill\break
\noindent {\footnotesize \bf malesevic@etf.bg.ac.yu}
\hfill}

\newpage

\end{document}